\documentstyle[12pt]{article}
\begin{document}

\vspace{1.cm}

\begin{center}
{\huge {\bf Isoperimetric Problems on cyclic polygons}} \footnote{
Mathematics Subject Classification 51M10, 51M25, 52A40} \\

\end{center}

\vspace{2.cm} {\large \ Abd Raouf Chouikha} \footnote{University of Paris-Nord\
LAGA,CNRS UMR 7539,\\ 
{\small chouikha@math.univ-paris13.fr}}

\vspace{3.cm}

{\bf Abstract}\newline
{\it In this paper we are interested in isoperimetric inequalities for plane n-gons inscribed in a circle. More specifically, we will examine a conjecture provided by Paul Levy that we propose to solve in particular cases.}\newline

\vspace{2.cm}

\section{Introduction}

In [4] we presented a problem formulated by Paul Levy concerning cyclic $n$-gons. He wanted to compare the area $A_n$ of a polygon $\Pi_n$ with the pseudo-area $P_n$ (here $L_n$ is its perimeter and $a_i$ are the length of sides)
\begin{equation}
P_n = {\frac{ L_n^2}{{4}}} \sqrt{ (1 - {\frac{2a_1}{{L_n}}})(1 - {\frac{2a_2%
}{{L_n}}})(1 - {\frac{2a_3}{{L_n}}}).....(1 - {\frac{2a_n}{{L_n}}})}
\label{6}
\end{equation}
It has been well known since antiquity that for $n=3$ or $n=4,$ we have perfect equality:
$$A_3 = P_3,\quad A_4 = P_4.$$
To that end Paul Levy defined the ratio $${\displaystyle \phi_n = {\frac{A_n}{{P_n}}}},$$
for any $n$-gon \ $\Pi_n$, with sides of lengths $a_1,a_2,....,a_n$
enclosing an area $A_n,$ and pseudo-area $P_n$. He remarked for the regular $n$-gon

\begin{equation}
\phi_n^0 = \frac{{A_n}^0}{{P_n}^0} = {\frac{1}{{{n \tan{\frac{\pi}{{n}}}}\ (1-{\frac{2}{{n}}}%
)^{n/2}}}} \label{7}
\end{equation}
He proposed the following\\

{\bf Conjecture (L):} {\it Define the ratio ${\displaystyle \phi_n = {\frac{
A_n}{{P_n}}}}$. For any $n$-gon \ $\Pi_n$, with sides of lengths $a_1,a_2,....,a_n$
enclosing an area $A_n,$ and $P_n$ defined as above, this ratio verifies }
$$\phi_n^0 < \phi_n < 1.$$\\
Several examples, initially proposed by Paul Levy himself and later improved and developed by [Ch], tend to prove the veracity of such a conjecture in a more general case.\\
In this paper we consider two types of approaches with the ultimate goal of resolving Conjecture (L).\\
 The first consists of deforming a regular polygon.
 At fist this conjecture is proved for a $n$-gon closed to the regular one. Subsequently, the following idea is explored: start with a regular $n$-gon with perimeter $L_n$, area $A_n$ and pseudo-perimeter $P_n$ and to deform each of its sides so that the new $n$-gon, with an invariant perimeter (equal to that of the regular one). More precisely, let \ $a= 2 R \sin \frac{\pi}{n}$\  be the length of the side,  each vertex of the regular $n$-gon undergoes a rotation of angle $y_i$ so that the new polygon has $n$ sides whose length is\
 $2 R \sin (\frac{\pi}{n}+y_i).$\ Moreover, it is cyclic and has the same perimeter as the regular one. We then prove for $-\frac{\pi}{4}<\frac{\pi}{n}+y_i <\frac{\pi}{4}, \ n > 4$ \ the left inequality of Conjecture (L) \  \quad $\phi_n^0 < \phi_n.$\\
 In the second approach, the previous method is applied to arbitrary cyclic polygons\ $\Pi$.\ Results can only be proven for a very slight deformation of \ $\Pi.$\ The more general case unfortunately still seems inaccessible.\\

\section{Deformation of regular $n$-gon}\quad 

Let us consider a regular  $n$-gon with sides of length\ $a = 2 R \sin\frac{\pi}{n}.$  Starting with the regular polygon, each vertex is rotated by an angle \ $y_i.$ \   This yields a polygon \ $\Pi_n$\ whose sides have length \ $a_i = 2 R \sin(\frac{\pi}{n}+y_i),$\ so that the perimeter remains unchanged \ $\sum _{i=1}^{n}\sin \left( {\frac {\pi }{n}}+y_{{i}} \right) =n \sin \left( {\frac {\pi }{n}} \right) .$\ The perimeter and area of \ $\Pi_n$\ are then $$L_n = 2 R \sum_{i=1}^n \sin(\frac{\pi}{n}+y_i)= 2 R n \sin \left( {\frac {\pi }{n}} \right), \ A_n = \frac{R^2}{2} \sum_{i=1}^n \sin(\frac{2\pi}{n}+2y_i).$$
The pseudo-area is
$$P_n={R}^{2} \left( \sum _{i=1}^{n}\sin \left( {\frac {\pi }{n}}+y_i
 \right)  \right) ^{2}\sqrt {\prod _{i=1}^{n} \left[  1- \frac{2\,\sin \left( 
{\frac {\pi }{n}}+y_i \right)}{  \left( \sum _{i=1}^{n}\sin \left( {
\frac {\pi }{n}}+y_i \right)  \right)} \right]  }.$$
Then the quotient is 
\begin{equation}\phi_n = \frac{A_n}{P_n} = \frac{1}{2}  \frac {\sum _{i=1}^{n}\sin \left( {\frac {2\pi }{n}}+2 y_i
 \right)}{ \left( \sum _{i=1}^{n}\sin \left( {\frac {\pi }{n}}+y_i
 \right)  \right) ^{2}\sqrt {\prod _{i=1}^{n} \left[ 1- \frac{2\,\sin \left( 
{\frac {\pi }{n}}+y_i \right)}{  \left( \sum _{i=1}^{n}\sin \left( {
\frac {\pi }{n}}+y_i \right)  \right)} \right] }}.\end{equation}\\

Let us write the length side of polygon \ $\Pi_n,$  $$a_i = \frac{L_n}{n} + x_i, \qquad x_i = 2 R \left(\sin(\frac{\pi}{n}+y_i) -\sin(\frac{\pi}{n})  \right).$$Since the perimeter remains unchanged then \ $\sum _{i=1}^{n} x_i\ =0.$\
Then, perimeter area and pseudo-area may be written 
$$L_n = 2 n R\sin \frac{\pi}{n}, \quad A_n =\frac{1}{4}\sum _{i=1}^{n} \left( {\frac {L_{{n}}}{n}}+x_{{i}} \right) 
\sqrt {4{R}^{2}- \left( {\frac {L_{{n}}}{n}}+x_{{i}} \right) ^{2}},$$
$$P_{{n}}=\frac{1}{4} {L_{{n}}}^{2}\sqrt {\prod _{i=1}^{n}\left[1-\frac{2}{n}-{ 
\frac {2 x_{{i}}}{L_{{n}}}}\right]}.$$ Then the  quotient is 
\begin{equation}\phi_n = {\frac{ \sum _{i=1}^{n}\left( \sin  \left({\frac {\pi }{n}}\right)
 +\frac{x_i}{2 R} \right) \sqrt {1-  \left(\,\sin \left( {
\frac {\pi }{n}}\right)  +\frac{x_i}{2 R} \right) ^{2}}}{{n}^{2}
 \left( \sin  {\frac {\pi }{n}}   \right) ^{2} 
\sqrt {\prod _{i=1}^{n}\left [1-\frac {2 }{n}-{\frac{x_{{i}}}{{n}\ {R}\  
\sin  {\frac {\pi }{n}}   }}\right]}}} .\end{equation}\\

{\bf Remarks 2.1}\\ 
{\bf 1 - }\ Formulas (3) or (4) seem to be the most general and the most suitable for proving Conjecture (L). As is readily apparent, it is difficult to use this formula as it stands to prove Conjecture (L).\\
{\bf 2 - }\ Instead of assuming the perimeter remains invariant, we could of course have assumed the area remains invariant. But the calculations seem little more complicated.

\subsection{Polygons closed to the regular one}\quad 
Let us examine the case where the new polygon is very close to the regular one. We then obtain the following result\\ 

{\bf Theorem 2.2} \quad {\it Let us consider a regular\ $n$-gon $ n\geq 4$\ inscribed in a circle of radius $R$ where each vertex is rotated by a small angle $y_i, i=1...n$ so that the perimeter remains unchanged. Then, the ratio $\phi_n$ of this new $n$-gon is such that 
$$\phi_n^0 < \phi_n =  \frac{1}{2}  \frac {\sum _{i=1}^{n}\sin \left( {\frac {2\pi }{n}}+2 y_i
 \right)}{ \left( \sum _{i=1}^{n}\sin \left( {\frac {\pi }{n}}+y_i
 \right)  \right) ^{2}\sqrt {\prod _{i=1}^{n} \left[ 1- \frac{2\,\sin \left( 
{\frac {\pi }{n}}+y_i \right)}{  \left( \sum _{i=1}^{n}\sin \left( {
\frac {\pi }{n}}+y_i \right)  \right)} \right] }} < 1 .$$
Conjecture (L) is then verified.}\\

{\bf Proof of Theorem 2.2}\\
From (4) the derivative of $\phi_n$ with respect to $x_i$ at $x_i=0$ is 
$${\frac{d \phi_n}{dx_i}} (0)  = \frac{n^{\frac{n-2}{2}}}{2 \left( n-2 \right) ^{\frac{n}{2}}} \frac{\left( 3\,{n} \left( 
\cos  {\frac {\pi }{n}}   \right) ^{2}-4
 \left( \cos  {\frac {\pi }{n}}   \right) ^{2}-{n}
+2 \right) }{ \left( n-2 \right)  \left( \cos
  {\frac {\pi }{n}}   \right) \left(  \left( 
\sin  {\frac {\pi }{n}}   \right) ^{2} \right){R}} > 0.$$
That means when the regular $n$-gon inscribed in a circle of radius $R$ is perturbed by slightly shifting its sides, the quotient of the resulting polygon \ $\Pi_n,$\ increases : \ $\phi_n^0 < \phi_n.$\ Moreover, since \ $x_i$\ is very small and \ $\phi_n^0 <1$. We then deduce \ $\phi_n < 1.$\ Thus, Conjecture (L) is verified.\\

\subsection{Other special cases of cyclic $n$-gon} \quad
However, Conjecture (L) can be proven for particular cases. For example, when each vertex is rotated by an angle \ $y_i =(-1)^i y$ \ from the regular $n$-gon.  This yields a cyclic polygon \ $\Pi_n$\ whose sides have length \ $a_i = 2 R \sin(\frac{\pi}{n}+(-1)^i y).$\ The perimeter and area of \ $\Pi_n$\ are then $$L_n = 2 R \sum_{i=1}^n \sin(\frac{\pi}{n}+(-1)^i y), \ A_n = \frac{R}{2} \sum_{i=1}^n \sin(\frac{2\pi}{n}+2(-1)^i y).$$
To go further, we must examine two cases separately.

\subsubsection{The number of sides is even} \quad In this case we obtain the perimeter and the area
$$L_n = 2 n R \sin\frac{\pi}{n} \cos y , \ A_n = \frac{1}{2} n R^2 \sin\frac{2\pi}{n} \cos 2y,$$
as well as the pseudo-area
$$P_n=n\tan {\frac {\pi }{n}}   \left( \cos y  \right) ^{2} \left( 1-\frac{4}{n}+\frac{4}{{n}^{2}}-\frac{4 \left( 
\tan  y    \right) ^{2}}{ \left( \tan  {\frac 
{\pi }{n}}   \right) ^{2}{n}^{2}} \right) ^{\frac{n}{4}}.$$
We then derive a simple form of the quotient
$$\phi_n = \frac{A_n}{P_n} = \frac{\cos \left(2y  \right)}{ {n} \left( \tan \left( {\frac {
\pi }{n}} \right)  \right)  \left( \cos \left( y  \right) 
 \right) ^{2}   \left( 1-\frac{4}{n}+\frac{4}{{n}^{2}}-\frac{4 \left( 
\tan  y    \right) ^{2}}{ \left( \tan  {\frac 
{\pi }{n}}   \right) ^{2}{n}^{2}} \right) ^{\frac{n}{4}} } .$$

 We easily deduce $$\phi_4 = {\frac {\cos \left( 2y \right) }{ 4\left( \cos \left( y \right) 
 \right) ^{2} \left( \frac{1}{4}-\frac{1}{4} \left( \tan  y   \right) ^
{2} \right) }} = 1. $$ For $y=0$ we have obviously \ $\phi_{n,0} = \phi_n.$\\ Moreover, by {\it Maple} we prove \ $\phi_n < 1$\ as well as
$$1 < \frac{\cos \left( {2y} \right)  \left( 1-\frac{2}{n} \right) ^{
\frac{n}{2}}}{ \left( \cos \left( {y} \right)  \right) ^{2}
   \left( 1-\frac{4}{n}+\frac{4}{{n}^{2}}-\frac{4 \left( 
\tan  y    \right) ^{2}}{ \left( \tan  {\frac 
{\pi }{n}}   \right) ^{2}{n}^{2}} \right) ^{\frac{n}{4}}},$$
which means \ $\phi_n^0 < \phi_n$. Thus, this polygon \ $\Pi_n$\ verifies conjecture (L).

\subsubsection{The number of sides is odd} \quad We consider here a $n$-gon \ $\Pi_n$\ with one side of length \ $a_n = R \sin\frac{\pi}{n}$.\ The others sides are such that \ $a_i= \sin(\frac{Pi}{n}+\frac{(-1)^i y}{n})$\ for \ $1 \leq i \leq n-1$.\ This yields
$$L_n = 2 R \sin\frac{\pi}{n} \left(1 + (n-1) \cos y  \right), \ A_n = \frac{1}{2} R^2 \sin\frac{2\pi}{n} \left(1 + (n-1) \cos 2y  \right).$$
We then compute the quotient
$$\phi_n =  \frac{ \left( 1+\cos 2y    \left( n-1 \right) 
 \right) }{ \left( \tan  {\frac {\pi }{n}}   \right)  
 \left( 1+ \left( n-1 \right) \cos  y    \right) ^
{2} \left[ \frac{   \left(  \left( 1+ \left( n-3 \right) \cos  
y    \right) ^{2}-4\, \left( \cot  {\frac {\pi }{n
}}   \right) ^{2} \left( \sin  y   
 \right) ^{2} \right)}{  \left( 1+ (n-1) \cos  y   
    \right) ^{2}} \right] ^{\frac{n-1}{4}}
 \sqrt {1-\frac{2}{ \left( 1+ \left( n-1 \right) 
\cos  y    \right)} }}.$$

For $y=0$ we have obviously \ $\phi_n^0 = \phi_n.$\\
We find of course for $n=3$
$$\phi_3
= \frac{ \left( 1+2\,\cos \left( 2\,y \right)  \right) \sqrt {3}}{3 \left( 1
+2\,\cos \left( y \right)  \right) ^{2}}{\frac {1}{\sqrt {{\frac {1-4/
3\, \left( \sin \left( y \right)  \right) ^{2}}{ \left( 1+2\,\cos
 \left( y \right)  \right) ^{2}}}}}}{\frac {1}{\sqrt {1-\frac{2}{ \left( 1+2
\,\cos \left( y \right)  \right)} }}} =1.$$
Moreover, by {\it Maple} it yields \ $\phi_n < 1$\ as well as
$$1 \leq \frac{n \left(1-\frac{2}{n}\right)^{\frac{n}{2}} \left( 1+\cos 2y    \left( n-1 \right) 
 \right) }{   
 \left( 1+ \left( n-1 \right) \cos  y    \right) ^
{2} \left[ \frac{   \left(  \left( 1+ \left( n-3 \right) \cos  
y    \right) ^{2}-4\, \left( \cot  {\frac {\pi }{n
}}   \right) ^{2} \left( \sin  y   
 \right) ^{2} \right)}{  \left( 1+ (n-1) \cos  y   
    \right) ^{2}} \right] ^{\frac{n-1}{4}}
 \sqrt {1-\frac{2}{ \left( 1+ \left( n-1 \right) 
\cos  y    \right)} }}.$$
This means \ $\phi_n^0 < \phi_n$.

We deduce from this that the conjecture (L) holds in this case as well. Thus, we proved:\\

{\bf Theorem 2.3} \quad {\it Let us consider a regular $n$-gon inscribed in a circle of radius $R$ for which each vertex is rotated by angle $(-1)^i y, i=1...n$ . Then, the obtained $n$-gon \ $\Pi_n$\ is cyclic with sides of length \ $a_i = 2 R \sin (\frac{\pi}{n}+(-1)^i y).$\ Moreover, \ $\Pi_n$\ satisfied Conjecture (L).}

\subsection{More general cases}
In fact, we may prove little more than above. Let us compute the derivative of $\phi_n$ with respect to $x_i$ (using {\it Maple})
$${\frac{d \phi_n}{dx_i}} = \frac{{L_{{n}}}^{\frac{n-4}{2}}{n}^{\frac{n-2}{2}}} { \left[ \prod _{i=1}^{n}\left( nL_{{n}}-2
\,L_{{n}}-2\,x_{{i}}n\right) \right] ^{\frac{3}{2}}} \ \psi_n,$$
where $$\psi_n =  2\,\sum _{i=1}^{n}{\frac {2
\,{R}^{2}{n}^{2}-{L_{{n}}}^{2}-2\,L_{{n}}x_{{i}}n-{x_{{i}}}^{2}{n}^{2}
}{\sqrt {4\,{R}^{2}{n}^{2}-{L_{{n}}}^{2}-2\,L_{{n}}x_{{i}}n-{x_{{i}}}^
{2}{n}^{2}}}}\prod _{i=1}^{n}\left(nL_{{n}}-2\,L_{{n}}-2\,x_{{i}}n\right)+$$ $$\sum _{i=
1}^{n} \left( L_{{n}}+x_{{i}}n \right) \sqrt {4\,{R}^{2}{n}^{2}-{L_{{n
}}}^{2}-2\,L_{{n}}x_{{i}}n-{x_{{i}}}^{2}{n}^{2}}\ \times $$ $$\sum _{i=1}^{n}
 \left( {\frac {\prod _{{\it j}=1}^{i} \left(   nL_{{n}}-2\,L_{
{n}}-2\,x_{{{\it j}}}n \right) \prod _{{\it j}=i}^{n} \left(nL_{{n}}-2\,L_{
{n}}-2\,x_{{{\it j}}}n \right) }{ \left( nL_{{n}}-2\,L_{{n}}-2\,x_{{i
}}n \right) ^{2}}} \right)  .$$

When examining expression of ${\frac{d \phi_n}{dx_i}}$ more closely, one notices that all the terms are positive except
$$2\,{R}^{2}{n}^{2}-{L_{{n}}}^{2}-2\,L_{{n}}x_{{i}}n-{x_{{i}}}^{2}{n}^{2},$$
this one is positive when for\ $1 \leq i\leq n$\ the following hold
$$- R\left(\sqrt2+2\sin \frac {\pi }{n}\right) < x_i < R\left(\sqrt2-2\sin \frac {\pi }{n}\right),$$
or equivalently since \ $x_i = 2 R \left(\sin(\frac{\pi}{n}+y_i) -\sin(\frac{\pi}{n})  \right)$
$$- \sqrt2 < 2 \sin(\frac{\pi}{n}+y_i) <  \sqrt2\quad \Leftrightarrow \quad -\frac{\pi}{4}<\frac{\pi}{n}+y_i <\frac{\pi}{4}.$$
We then proved\\

{\bf Theorem 2.4}\quad {\it  Let us consider a regular  $n$-gon $n > 4,$ inscribed in a circle of radius $R$ with sides of length\ $a = 2 R \sin\frac{\pi}{n}.$ The inscribed polygon $\Pi_n$ obtained from the regular one, in which each vertex undergoes a rotation by an angle \ $y_i$ \ and the perimeter remains unchanged is such that :\\
when \quad $-\frac{\pi}{4}<\frac{\pi}{n}+y_i <\frac{\pi}{4}$\quad then \quad $\phi_n^0 < \phi_n.$}\\

\section{Deformation of a cyclic $n$-gon}\quad
This section can be seen as an improvement on the preceding section, which only considered a regular $n$-gons.\\
Let us consider a cyclic non regular $n$-gon with sides of length\ $a_i = 2 R \sin\theta_i.$  Starting with this polygon, each vertex is rotated by an angle \ $y_i.$ \   This yields a cyclic polygon \ $\Pi_n$\ whose sides have length \ $a_i = 2 R \sin(\theta_i+y_i),$\ so that the perimeter remains unchanged \ $\sum _{i=1}^{n} \sin \theta_i = \sum _{i=1}^{n}\sin (\theta_i+y_i) .$\ Let us write the length side of polygon \ $\Pi_n,$  $$b_i = a_i + x_i, \qquad x_i = 2 R \left(\sin(\theta-i+y_i) -\sin \theta_i  \right).$$Since the perimeter remains unchanged then \ $\sum _{i=1}^{n} x_i\ =0.$\
Then, perimeter area and pseudo-area may be written 
$$L_n = \sum _{i=1}^{n} b_i = \sum _{i=1}^{n} a_i, \quad A_n = \sum _{i=1}^{n}\frac{1}{4} \left( a_{{i}}+x_{{i}} \right) \sqrt {4\,
{R}^{2}-{a_{{i}}}^{2}-2\,a_{{i}}x_{{i}}-{x_{{i}}}^{2}},$$
$$P_{{n}}=\frac{1}{4} \left( \sum _{i=1}^{n}a_{i} \right) ^{2}\sqrt 
{\prod _{i=1}^{n} \left( 1-{\frac {2\,a_{{i}}+2\,x_{{i}}}{\sum _{i=1}^
{n}a_{{i}}}} \right) }.$$
Then the quotient is 
\begin{equation}\phi_n = {\frac{\sum _{i=1}^{n} \left( a_{{i}}+x_{{i}} \right) \sqrt {4\,
{R}^{2}-{a_{{i}}}^{2}-2\,a_{{i}} x_{{i}}-{x_{{i}}}^{2}}}{\left[ \sum _{i=1}^{n}a_{{i}} \right] ^{2}\sqrt 
{\prod _{i=1}^{n} \left( 1-{\frac {2\,a_{{i}}+2\,x_{{i}}}{\sum _{i=1}^
{n}a_{{i}}}} \right) }}}\end{equation}

\subsection{A special case}\quad Let us consider now a cyclic $n$-gon with $n$ even sides equal in pairs :\ $a_i = a_{i-1+\frac{n}{2}}, i=1... \frac{n}{2},$\ ($n$ sufficiently large).  We deform each side so that the new polygon is always cyclic, with side lengths of $a_i+x_i=a_i+(-1)^i x$ where $x$ is small. Then the perimeter remains unchanged since \ $\sum _{i=1}^{n} x_i\ =0:$
 $$\ L_n = \sum _{i=1}^{n}a_{{i}}+ \left( -1 \right) ^{i} x = \sum _{i=1}^{n}a_{{i}} = 2 \sum _{i=1}^{\frac{n}{2}}a_{{i}}.$$
Area and pseudo-area may be written 
$$ A_n = \frac{1}{4}\sum _{i=1}^{\frac{n}{2}} \left( a_{{i}}+x \right) \sqrt {4\,{R}^{2}-{a_{{i
}}}^{2}-2\,a_{{i}}x-{x}^{2}}+\frac{1}{4}\sum _{i=1}^{\frac{n}{2}} \left( a_{{i}}-x
 \right) \sqrt {4\,{R}^{2}-{a_{{i}}}^{2}+2\,a_{{i}}x-{x}^{2}},$$
$$P_{{n}}=\frac{1}{4}{L_{{n}}}^{2}\sqrt {\prod _{i=1}^{\frac{n}{2}} \left(1-{\frac {4 a_{
{i}}}{L_{{n}}}}+{\frac {4{a_{{i}}}^{2}}{{L_{{n}}}^{2}}}-{\frac {4{
x}^{2}}{{L_{{n}}}^{2}}}\right)}.$$
The quotient \ $\frac{A_n}{P_n}$\ is 
$$\phi_n =  \frac{ \sum _{i=1}^{\frac{n}{2}}\left[ \left( a_{{i}}+x \right) \sqrt {4\,{R}^{2
}-{a_{{i}}}^{2}-2\,a_{{i}}x-{x}^{2}}+ \left( a_{{i
}}-x \right) \sqrt {4\,{R}^{2}-{a_{{i}}}^{2}+2\,a_{{i}}x-{x}^{2}}
  \right]}{{L_{{n}}}^{2} \sqrt {\prod _{i=1}^{\frac{n}{2}} \left(1-{
\frac {4 a_{{i}}}{L_{{n}}}}+{\frac {4{a_{{i}}}^{2}}{{L_{{n}}}^{2}}}-{\frac {4{x}^{2}}{{L_{{n}}}^{2}}}\right)}}.$$

We then propose\\

{\bf Theorem 2.5}\quad {\it  Let us consider a non regular polygon denoted \ $\Pi_{n,0}$\ with $n$ even sides equal in pairs and inscribed in a circle of radius $R$.  Consider the inscribed polygon $\Pi_{n,x}$ obtained from this one, in which each vertex undergoes a rotation by a small angle \ $y_i$ \ so that the perimeter \ $L_n$\ remains unchanged with sides of length $a_i + (-1)^i x, i=1... \frac{n}{2}, $ \  $x$ very small.\\
Then for $n$ large we have \ $ \phi_{n,0}=\frac{A_{n,0}}{P_{n,0}} < \frac{A_{n,x}}{P_{n,x}}$\
where \ $\phi_{n,0}$ and \ $\phi_{n,x}$\  are the quotients of \ $\Pi_{n,0}$\ and \ $\Pi_{n,x}$\ respectively.\\ 
If in addition \ $\phi_{n,0} < 1$\ then for small\ $x$ \ $\phi_{n,x} < 1$.}\\

This means, according to this theorem, that if \ $\Pi_{n,0}$\ verified Conjecture (L) so the same applies to \ $\Pi_{n,x}.$\\

{\bf Proof of Theorem 2.5} \quad Consider the series expansion to order 3 of \ $\phi_n$\ with respect to $x$ in the neighborhood of $0$
$$\phi_n = \phi_{n,0} + \phi_{n,1} x + \phi_{n,2} x^2 + O(x^3),$$
where \ $\phi_{n,0}$\ is the quotient of \ $\Pi_{n,0}$\ and
$$\phi_{n,1} =  \frac{\left( \sum _{i=1}^{n}-{\frac {{a_{{i}}}^{2}}{\sqrt {4\,{R}^{2}-{a_{{
i}}}^{2}}}}+\sqrt {4\,{R}^{2}-{a_{{i}}}^{2}}+\sum _{i=1}^{n}{\frac {{a
_{{i}}}^{2}}{\sqrt {4\,{R}^{2}-{a_{{i}}}^{2}}}}-\sqrt {4\,{R}^{2}-{a_{
{i}}}^{2}} \right)}{ {L_n^2 \sqrt{\prod _{i=1
}^{\frac{n}{2}} \left(1-4\,{\frac {a_{{i}}}{L_{{n}}}}+4\,{\frac {{a_{{i}}}^{2}}{{L
_{{n}}}^{2}}} \right)}}
} =0.$$
$$\phi_{n,2} = -\frac{\sum _{i=1}^{\frac{n}{2}}a_{{i}}\sqrt {4\,{R}^{2}-{a_{{i}}}^{2}}\ \sum _{i=1
}^{\frac{n}{2}}\frac{-4}{ \left( {L_{{n}}}^{2}-4\,a_{{i}}L_{{n}}+4\,{a_{{i}}}^{2}
 \right)}}{ L_n^2 \sqrt{\prod _{i=1}^{\frac{n}{2}} \left(1-4\,{\frac {a_{{i}}}{L_{{n
}}}}+4\,{\frac {{a_{{i}}}^{2}}{{L_{{n}}}^{2}}} \right)}} +
$$ $$\frac{2\,\sum _{i=1}^{\frac{n}{2}}a_{{i}}\sqrt {4\,{R}^{2}-{a_{{i}}}^{2}} \left( 
- \frac{1}{\left( 8\,{R}^{2}-2\,{a_{{i}}}^{2} \right)}-{\frac {{a_{{i
}}}^{2}}{ \left( 8\,{R}^{2}-2{a_{{i}}}^{2} \right) ^{2}}} \right) -{
\frac {a_{{i}}}{\sqrt {4\,{R}^{2}-{a_{{i}}}^{2}}}}}{L_n^2 \sqrt{ \prod _{i=1
}^{\frac{n}{2}} \left(1-4\,{\frac {a_{{i}}}{L_{{n}}}}+4\,{\frac {{a_{{i}}}^{2}}{{L
_{{n}}}^{2}}} \right)} } = \frac{N}{D}.$$

The numerator $N$ of \ $\phi_{n,2}$\ can be expressed in a simpler way by the following\\
                 
{\bf Lemma 2.6}\ {\it The numerator $N$ can be written 
$$N=  -4\,\sum _{i=1}^{n}a_{{i}}\sqrt {4\,{R}^{2}-{a_{{i}}}^{2}}\ 
\sum _{i=1}^{\frac{n}{2}} \frac{1}{\left( L_{{n}}-2\,a_{{i}} \right) ^{2}}+2\,\sum _
{i=1}^{n}{\frac { \left( 6\,{R}^{2}-{a_{{i}}}^{2} \right) a_{{i}}}{
 \left( 4\,{R}^{2}-{a_{{i}}}^{2} \right) ^{\frac{3}{2}}}}.$$
Moreover, the following inequalities hold}  $$1-6\,\sum _{i=1}^{\frac{n}{2}} \frac{1}{ \left( L_{{n}}-2\,a_{{i}} \right) ^{2}} >0, \quad 
 \sum _{i=1}^{n} 3 a_i{\frac {6\,{R}^{2}-{a_{{i}}}^{2}}{ \left( 4\,{R}^{2}-{a_{{i}}}^{2}
 \right) ^{\frac{3}{2}}}}-\sum _{i=1}^{n} a_i\sqrt {4\,{R}^{2}-{a_{{i}}}^{2}} > 0,$$
and they imply \ $N > 0.$\\

{\bf Proof of Lemma 2.6} \quad let us prove at first for \ $0 < a_i < 2R$\ that $$ 3 a_i{\frac {6\,{R}^{2}-{a_{{i}}}^{2}}{ \left( 4\,{R}^{2}-{a_{{i}}}^{2} \right) ^{\frac{3}{2}}}}- a_i\sqrt {4\,{R}^{2}-{a_{{i}}}^{2}} > 0.$$
Indeed, that expression can be written $$\frac{3 (6 R^2 - a_i^2) - (4 R^2 - a_i^2)^2}{\sqrt {4\,{R}^{2}-{a_{{i}}}^{2}}}.$$
So, it easy to see that \ $6 R^2 + 3 y - y^2$\ is positive where \\ $y = 4 R^2 - a_i^2, 0 < y < 4R.$\\
Furthermore, since \ $\Pi_{n,0}$\ is not regular, then it has at least one side of length \ $ a_p$ \ less than \ $ a_i, i=1...n$\  and at least one side of length \ $ a_q$ \ greater than \ $ a_i, i=1...n$. So, \ $ a_p < a_i < a_q$.\\
We then deduce \ $L_n -2 a_i > L_n -2 a_q$\ which implies 
$$1-3\,{\frac {n}{ \left( L_{{n}}-2\,a_{{q}} \right) ^{2}}}<1-6\,\sum _{i=
1}^{\frac{n}{2}} \frac{1}{4\left( L_{{n}}-2\,a_{{i}} \right) ^{2}}.$$
Moreover, for $n$ sufficiently large the left part of that inequality is positive \ $0 < 1-3\,{\frac {n}{ \left( L_{{n}}-2\,a_{{q}} \right) ^{2}}}.$\\ Thus, 
$$N > 4\,\sum _{i=1}^{n}a_{{i}}\sqrt {4\,{R}^{2}-{a_{{i}}}^{2}} \left(1-6\,\sum _{i=
1}^{\frac{n}{2}} \frac{1}{4\left( L_{{n}}-2\,a_{{i}} \right) ^{2}}\right) > 0.$$

\end{document}